\documentclass[leqno]{amsart}
\usepackage{amsmath}
\usepackage{amssymb}
\usepackage{amsthm}
\usepackage{enumerate}
\usepackage[mathscr]{eucal}
\theoremstyle{plain}
\newtheorem{theorem}{Theorem}[section]
\newtheorem{prop}[theorem]{Proposition}
\newtheorem{corollary}{Corollary}[theorem]

\newtheorem{conj}[theorem]{Conjecture}

\newtheorem{remark}[theorem]{Remark}
\theoremstyle{definition}

\theoremstyle{remark}
\newtheorem{example}{Example}

\numberwithin{equation}{section}
\begin{document}

\title[Birkhoff-James orthogonality of linear operators on finite dimensional Banach spaces]{Birkhoff-James orthogonality of linear operators on finite dimensional Banach spaces}
\author[ Debmalya Sain  ]{ Debmalya Sain  }

\newcommand{\acr}{\newline\indent}

\address{\llap{\,}Department of Mathematics\acr
                              Jadavpur University\acr
                              Kolkata 700032\acr
                              West Bengal\acr
                              INDIA}
\email{saindebmalya@gmail.com}

\thanks{} 


\begin{abstract}
In this paper we characterize Birkhoff-James orthogonality of linear operators defined on a finite dimensional real Banach space $ \mathbb{X}. $ We also explore the symmetry of Birkhoff-James orthogonality of linear operators defined on $ \mathbb{X}. $ Using some of the related results proved in this paper, we finally prove that $ T \in \mathbb{L}(l_{p}^2) (p \geq 2, p \neq \infty) $ is left symmetric with respect to Birkhoff-James orthogonality if and only if $ T $ is the zero operator. We conjecture that the result holds for any finite dimensional strictly convex and smooth real Banach space $ \mathbb{X}, $ in particular for the Banach spaces $ l_{p}^{n} (p > 1, p \neq \infty). $
\end{abstract}

\maketitle
\section{Introduction.} 

Birkhoff-James orthogonality \cite{2} plays a vital role in the study of geometry of Banach spaces. One of the prominent reasons behind this is the natural connection shared by Birkhoff-James orthogonality with various geometric properties of the space, like smoothness, strict convexity etc. Recently in \cite{5}, Sain and Paul have characterized finite dimensional real Hilbert spaces among finite dimensional real Banach spaces in terms of operator norm attainment, using the notion of Birkhoff-James orthogonality. More recently, symmetry of Birkhoff-James orthogonality of linear operators defined on a finite dimensional real Hilbert space $ \mathbb{H} $ has been explored by Ghosh et al. in \cite{3}. However, it was remarked in \cite{3} that analogous results corresponding to the far more general setting of Banach spaces remain unknown. The aim of the present paper is twofold: we characterize Birkhoff-James orthogonality of linear operators defined on a finite dimensional real Banach space $ \mathbb{X} $ and we also explore the symmetry of Birkhoff-James orthogonality of linear operators defined on $ \mathbb{X} $. Using some of the results proved in this paper, we finally study the ``left symmetry" of Birkhoff-James orthogonality of linear operators defined on $ l_{p}^{2} (p \geq 2, p \neq \infty). $

\noindent Let $ (\mathbb{X},\|\|) $ be a finite dimensional real Banach space. Let $ B_\mathbb{X}=\{x \in \mathbb{X} : \|x\| \leq 1\} $ and $ S_\mathbb{X}=\{x \in \mathbb{X} : \|x\|=1\} $ be the unit ball and the unit sphere of the Banach space $ \mathbb{X}$ respectively. Let $ \mathbb{L}(\mathbb{X}) $ denote the Banach space of all linear operators on $ \mathbb{X}, $ endowed with the usual operator norm.\\
\noindent For any two elements $ x,y \in \mathbb{X}, $ $ x $ is said to be orthogonal to $ y $ in the sense of Birkhoff-James, written as $ x \perp_B y, $ if 
\[ \|x\| \leq \|x+\lambda y\| ~ \mbox{for all}~ \lambda \in \mathbb R. \]
  Likewise, for any two elements $ T,A \in  \mathbb{L}(\mathbb{X}), $ $ T $ is said to be orthogonal to $ A $ in the sense of Birkhoff-James, written as $ T \perp_B A, $ if 
  \[ \|T\| \leq \|T+ \lambda A\| ~\mbox{for all}~ \lambda \in \mathbb R. \] 
In a finite dimensional Hilbert space $ \mathbb{H}, $ Bhatia and Semrl \cite{1} proved that for any two elements $ T, A \in \mathbb{L}(\mathbb{H}), $ $ T \perp_B A $ if and only if there exists $ x \in \mathbb{H} $ with $ \|x\|=1 $ such that $ \|Tx\|=\|T\| $ and $ Tx \perp_B Ax. $ Sain and Paul \cite{5} generalized the result for linear operators defined on  finite dimensional real Banach spaces in Theorem $ 2.1 $ of \cite{5} by proving the following result: \\ 

\noindent Let $ \mathbb{X} $ be a finite dimensional real Banach space. Let $T \in \mathbb{L}(\mathbb{X})$ be such that $ T $ attains norm only at $ \pm D, $ where $ D $ is a closed, connected subset of $ S_\mathbb{X}. $ Then for $A \in \mathbb{L}(\mathbb{X})$ with $ T \perp_B A, $  there exists $x \in D$ such that $ Tx \perp_B Ax $. \\

\noindent However, useful as it indeed is, the above result does not completely characterize Birkhoff-James orthogonality of linear operators defined on a finite dimensional real Banach space $ \mathbb{X}. $ For a linear operator $ T $ defined on a Banach space $ \mathbb{X}, $ let $ M_T $ denote the collection of all unit vectors in $ \mathbb{X} $ at which $ T $ attains norm, i.e., 

\[ M_T = \{  x \in S_{\mathbb{X}} : \| Tx \| = \| T \|  \}. \]

\noindent Since there exists a linear operator $ T \in \mathbb{L}(\mathbb{X}) $ such that $ M_T $ is not of the form $ \pm D,  $ where $ D $ is a closed, connected subset of $ S_\mathbb{X}, $ the above result can not be effectively applied, at least in its original form, to determine whether $ T \perp_B A $ for some $ A \in \mathbb{L}(\mathbb{X}). $ The following example, given in \cite{6}, validates our claim.

\begin{example}
Consider $(\mathbb {R}^2,\|.\|)$ whose unit sphere is the
regular hexagon with vertices at $\pm (1,0)$,
$\pm (\frac{1}{2},\frac{\sqrt{3}}{2})$,
$\pm (-\frac{1}{2},\frac{\sqrt{3}}{2})$.
Let
\[
T =
 	\begin{pmatrix}
  		\frac{3}{4} & -\frac{\sqrt{3}}{4} \\ 
  		\rule{0em}{3ex}\frac{\sqrt{3}}{4} & \hphantom{-}\frac{3}{4}
 	\end{pmatrix}.
\]
Then $\|T\|=1$ and $T$ attains norm only at the points $\pm (1,0)$,
$\pm (\frac{1}{2},\frac{\sqrt{3}}{2})$, $\pm (-\frac{1}{2},\frac{\sqrt{3}}{2})$ and hence $ M_T $ can not be of the form $ \pm D, $ where $ D $ is a closed, connected subset of $ S_\mathbb{X}. $
\end{example}

\noindent In this paper we introduce a particular notion, motivated by geometric observations made by us, in order to completely characterize Birkhoff-James orthogonality of linear operators defined on finite dimensional real Banach spaces. For any two elements $ x, y $ in a real normed linear space $ \mathbb{X}, $ let us say that $ y \in x^{+} $ if $ \| x + \lambda y \| \geq \| x \| $ for all $ \lambda \geq 0. $ Accordingly, we say that $ y \in x^{-} $ if $ \| x + \lambda y \| \geq \| x \| $ for all $ \lambda \leq 0. $ Using this notion we completely characterize Birkhoff-James orthogonality of linear operators defined on finite dimensional real Banach spaces.\\
\\
\noindent Next we consider the symmetry of Birkhoff-James orthogonality of linear operators defined on a finite dimensional real Banach space $ \mathbb{X}. $ For an element $ x \in  \mathbb{X}, $ let us say that $ x $ is left symmetric (with respect to Birkhoff-James orthogonality) if $ x \perp_{B} y $ implies $ y \perp_{B} x $ for any $ y \in \mathbb{X}. $ Similarly, let us say that $ x $ is right symmetric (with respect to Birkhoff-James orthogonality) if $ y \perp_{B} x $ implies $ x \perp_{B} y $ for any $ y \in \mathbb{X}. $ It was proved in \cite{3} that if $ \mathbb{H} $ is a finite dimensional real Hilbert space then $ T \in \mathbb{L}(\mathbb{H}) $ is a left symmetric point if and only if $ T $ is the zero operator and $ T \in \mathbb{L}(\mathbb{H}) $ is a right symmetric point if and only if $ M_T = S_{\mathbb{H}} $. In this paper we consider the problem in the more general setting of real Banach spaces and prove some related results corresponding to the symmetry of linear operators defined on a finite dimensional real Banach space. We give example to show that in a finite dimensional real Banach space $ \mathbb{X} $, which is not  a  Hilbert space, there may exist nonzero linear operators $ T \in \mathbb{L}(\mathbb{X}) $ such that $ T $ is a left symmetric point in $ \mathbb{L}(\mathbb{X}). $ Finally, using some of the results proved in this paper, we prove that $ T \in \mathbb{L}(l_{p}^{2}) (p \geq 2, p \neq \infty) $ is left symmetric if and only if $ T $ is the zero operator. Motivated by this result, we conjecture that for any finite dimensional strictly convex and smooth real Banach space $ \mathbb{X}, $ $ T \in \mathbb{L}(\mathbb{X}) $ is left symmetric if and only if $ T $ is the zero operator.

\section{Main results.}
\noindent In order to characterize Birkhoff-James orthogonality of linear operators defined on finite dimensional real Banach spaces, we have introduced the notions $ y \in x^{+}  $ and $ y \in x^{-},  $ for any two elements $ x, y $ in a real normed linear space $ \mathbb{X}. $ First we state some obvious but useful properties of this notion which would be used later on in this paper, without giving explicit proofs.\\

\begin{prop} Let $ \mathbb{X}  $ be a real normed linear space and $ x, y \in \mathbb{X}. $ Then the following are true: \\
(i) Either $ y \in x^{+}  $ or $ y \in x^{-}. $\\
(ii) $ x \perp_{B} y $ if and only if $ y \in x^{+}  $ and $ y \in x^{-}. $ \\
(iii) $ y \in x^{+} $ implies that $ \eta y \in (\mu x)^{+} $ for all $ \eta , \mu >0. $ \\
(iv) $ y \in x^{+} $ implies that $ -y \in x^{-} $ and $ y \in (-x)^{-}. $\\
(v) $ y \in x^{-} $ implies that $ \eta y \in (\mu x)^{-} $ for all $ \eta , \mu >0. $ \\
(vi) $ y \in x^{-} $ implies that $ -y \in x^{+} $ and $ y \in (-x)^{+}. $
\end{prop}

\noindent In the next theorem we use this notion to give a characterization of Birkhoff-James orthogonality of linear operators defined on a finite dimensional real Banach space.\\ 

\begin{theorem} Let $ \mathbb{X} $ be a finite dimensional real Banach space. Let $T, A \in \mathbb{L}(\mathbb{X}).$ Then  $T \perp_B A $  if and only if there exists $x, y \in M_T $ such that $ Ax \in Tx^{+} $ and $ Ay \in Ty^{-} $.
\end{theorem}
\noindent \textbf{Proof :} Let us first prove the easier ``if" part.\\
 Suppose there exists $x, y \in M_T $ such that $ Ax \in Tx^{+} $ and $ Ay \in Ty^{-}. $ For any $ \lambda \geq 0, \| T + \lambda A \| \geq \| Tx + \lambda Ax \| \geq \| Tx \| = \| T \|. $ Similarly, for any $ \lambda \leq 0, \| T + \lambda A \| \geq \| Ty + \lambda Ay \| \geq \| Ty \| = \| T \|. $ This proves that $T \perp_B A. $\\
\noindent Let us now prove the comparatively trickier ``only if" part.\\
\noindent Let $T, A \in \mathbb{L}(\mathbb{X})$ be such that $ T \perp_B A. $ If possible suppose that there does not exist $ x, y \in M_T $ such that $ Ax \in Tx^{+} $ and $ Ay \in Ty^{-} $. Using $ (i) $ of Proposition $ 2.1 $, it is easy to show that either of the following is true: \\
\noindent (i) $ Ax \in Tx^{+} $ for each $ x \in M_T $ and $ Ax \notin Tx^{-} $ for any $ x \in M_T $ \\
\noindent (ii) $ Ax \in Tx^{-} $ for each $ x \in M_T $ and $ Ax \notin Tx^{+} $ for any $ x \in M_T. $ \\
\noindent Without loss of generality, let us assume that $ Ax \in Tx^{+} $ for each $ x \in M_T $ and $ Ax \notin Tx^{-} $ for any $ x \in M_T. $ Consider the function $ g : S_\mathbb{X} \times [ -1 , 1] \longrightarrow  \mathbb{R} $ defined by 
$$g(x,\lambda) = \|Tx + \lambda Ax\|.$$
It is easy to check that $ g $ is continuous. Given any $ x \in M_T, $ since $ Ax \notin Tx^{-}, $ there exists $ \lambda_{x} < 0 $ such that $ g(x,\lambda_x) = \|Tx + \lambda_x Ax\| < \| Tx \| = \| T \|. $ By continuity of $ g, $ there exists $r_x , \delta_x > 0$ such that 
 $$ g(y, \lambda) < \|T\| ~~~ \mbox{for all} ~~~ y \in B(x, r_x) \cap S_\mathbb{X} ~\mbox{ and}~ \mbox{for all} ~ \lambda \in ( \lambda_x - \delta_x ,  \lambda_x + \delta_x).$$
Using the convexity property of the norm function, it is easy to show that $ g(y,\lambda) = \|Ty + \lambda Ay\| < \|T\| ~~~ \mbox{for all} ~~~ y \in B(x, r_x) \cap S_\mathbb{X} $ and $ \mbox{for all}~ \lambda \in (\lambda_x, 0).$ \\

\noindent  For any $ z \in S_\mathbb{X} \diagdown M_T, $  we have $ g(z,0)=\|Tz\|<\|T\|. $ Thus by continuity of $ g, $ there exist $ r_z, \delta_z > 0 $ such that $ g(y,\lambda)=\|Ty+\lambda Ay\|<\|T\| ~ \mbox{for all} ~ y \in B(z,r_z)\cap S_\mathbb{X} $ and $ \mbox{for all}~ \lambda \in ( - \delta_z , \delta_z). $ \\

\noindent Consider the open cover $ \{ B(x, r_x) \cap S_\mathbb{X} : x \in M_T\} \cup \{ B(z, r_z) \cap S_\mathbb{X} : z \in S_\mathbb{X} \diagdown M_T \} $ of $ S_\mathbb{X}.$ Since $ \mathbb{X}  $ is finite dimensional, $ S_\mathbb{X} $ is compact. This proves that the considered open cover has a finite subcover and so we get, 
  \[ S_X \subset \cup_{i=1}^{n_1} B(x_i,r_{x_i})  \cup_{k=1}^{n_2} B(z_k,r_{z_k}) \cap S_{\mathbb{X}},\]
  for some positive integers $  n_1, n_2, $ where each $ x_i \in M_T $ and each $ z_k \in S_\mathbb{X} \diagdown M_T. $\\
	
	\noindent Choose $ \lambda_0 \in (\cap_{i=1}^{n_1} (\lambda_{x_{i}}, 0) )~ \bigcap ~(\cap_{k=1}^{n_2} (-\delta_{z_{k}},\delta_{z_{k}})). $\\
\noindent Since $ \mathbb{X} $ is finite dimensional, $ T + \lambda_0 A $ attains its norm at some $ w_0 \in S_\mathbb{X} $. Either $ w_0 \in B(x_i,r_{x_i}) $ for some $ x_i \in M_T $ or $ w_0 \in B(z_k,r_{z_k}) $ for some $ z_k \in S_\mathbb{X} \diagdown M_T. $ In either case, it follows from the choice of $ \lambda_0 $ that $ \|T + \lambda_0 A \| = \|(T + \lambda_0 A) w_0 \| < \|T\|, $ which contradicts our primary assumption that $ T \perp_B A $ and thereby completes the proof of the ``only if" part. \\

\noindent Theorem $ 2.1 $ of Sain and Paul \cite{5} can be deduced as a corollary to the previous theorem.

\begin{corollary}
Let $ \mathbb{X} $ be a finite dimensional real Banach space. Let $T \in \mathbb{L}(\mathbb{X})$ be such that $ T $ attains norm only at $ \pm D, $ where $ D $ is a closed, connected subset of $ S_\mathbb{X}. $ Then for $ A \in \mathbb{L}(\mathbb{X}) $ with $ T \perp_B A, $  there exists $x \in D$ such that $ Tx \perp_B Ax $.
\end{corollary}

\noindent \textbf{Proof :} Since $ T \perp_B A, $ applying Theorem $ 2.2, $ we see that there exists $ x, y \in M_T = \pm D $ such that $ Ax \in Tx^{+} $ and $ Ay \in Ty^{-}. $ Moreover, it is easy to see that by applying $ (iv) $ and $ (vi) $ of Proposition $ 2.1, $ we may assume without loss of generality that $ x, y \in D. $ Then following the same line of arguments, as in Theorem $ 2.1 $ of \cite{5}, it can be proved that there exists $ u_0 \in D $ such that $ Au_0 \in Tu_{0}^{+} $ and $ Au_0 \in Tu_{0}^{-}, $ by using the connectedness of $ D $. However, this is equivalent to $ Tu_0 \perp_{B} Au_0, $ completing the proof of Theorem $ 2.1 $ of \cite{5}.

\begin{remark}
The previous theorem gives a complete characterization of Birkhoff-James orthogonality of linear operators defined on a finite dimensional real Banach space $ \mathbb{X} $. Moreover, as we will see later on in this paper, the theorem is very useful computationally as well as from theoretical point of view. It should be noted that the main idea of the proof of Theorem $ 2.2 $ was already there in the proof of Theorem $ 2.1 $ of \cite{5}. However, complete characterization of Birkhoff-James orthogonality of linear operators on $ \mathbb{X} $ could not be obtained in \cite{5}. This reveals the usefulness of the notion introduced by us in this paper to meet this end.
\end{remark}

\noindent Next we consider the symmetry of Birkhoff-James orthogonality of linear operators defined on a finite dimensional real Banach space $ \mathbb{X}. $ $ T \in \mathbb{L}(\mathbb{X}) $ is left symmetric if $ T \perp_B A $ implies that $ A \perp_B T $ for any $ A \in \mathbb{L}(\mathbb{X}). $ Similarly, $ T \in \mathbb{L}(\mathbb{X}) $ is right symmetric if $ A \perp_B T $ implies that $ T \perp_B A $ for any $ A \in \mathbb{L}(\mathbb{X}). $ In the following theorem we establish a useful connection between left symmetry of an operator $ T \in \mathbb{L}(\mathbb{X}) $ and left symmetry of points in the corresponding norm attainment set $ M_T. $ \\

\begin{theorem} Let $ \mathbb{X} $ be a finite dimensional strictly convex real Banach space. If $ T \in \mathbb{L}(\mathbb{X})$ is a left symmetric point then for each $ x \in M_T, Tx $ is a left symmetric point. 
\end{theorem}

\noindent \textbf{Proof :} First we observe that the theorem is trivially true if $ T $ is the zero operator. Let $ T $ be nonzero.~Since $ \mathbb{X} $ is finite dimensional, $ M_T $ is nonempty. If possible suppose that there exists $ x_1 \in M_T $ such that $ Tx_1 $ is not a left symmetric point. Since $ T $ is nonzero, $ Tx_1 \neq 0.$     Then there exists $ y_1 \in S_\mathbb{X} $ such that $ Tx_1 \perp_B y_1 $ but $ y_1 \not\perp_B Tx_1. $ Since $ \mathbb{X} $ is strictly convex, $ x_1 $ is an exposed point of the unit ball $ B_{\mathbb{X}}. $ Let $ H $ be the hyperplane of codimension $ 1 $ in $ \mathbb{X} $ such that $ x_1 \perp_B H. $ Clearly, any element $ x $ of $ \mathbb{X} $ can be uniquely written in the form $ x = \alpha_1 x_1 + h, $ where $ \alpha_1 \in \mathbb{R} $ and $ h \in H. $ Define a linear operator $ A \in \mathbb{L}(\mathbb{X}) $ as follows: 
\[ Ax_1 = y_1, Ah = 0 ~\mbox{for all}~ h \in H. \]

\noindent Since $ x_1 \in M_T $ and $ Tx_1 \perp_B Ax_1, $ it follows that $ T \perp_B A. $ Since $ T $ is left symmetric, it follows that $ A \perp_B T. $\\
\noindent It is easy to check that $ M_A = \{ \pm x_1 \}, $ since $ \mathbb{X} $ is strictly convex. Applying Theorem $ 2.1 $ of \cite{5} to $ A, $ it follows from $ A \perp_B T $ that $ Ax_1 \perp_B Tx_1, $ i.e., $ y_1 \perp_B Tx_1, $ contrary to our initial assumption that $ y_1 \not\perp_B Tx_1. $ This contradiction completes the proof of the theorem. \\

\begin{corollary}
Let $ \mathbb{X} $ be a finite dimensional real Banach space such that the unit sphere $ S_{\mathbb{X}} $ has no left symmetric point. Then $ \mathbb{L}(\mathbb{X}) $ can not have any nonzero left symmetric point. In particular, if $ \mathbb{H} $ is a finite dimensional real Hilbert space then $ \mathbb{L}(\mathbb{L}(\mathbb{H})) $ has no nonzero left symmetric point.
\end{corollary}

\noindent \textbf{Proof :} The first statement follows from the previous theorem and the fact that every linear operator defined on a finite dimensional Banach space must attain norm at some point of the unit sphere. The second statement follows from the first statement and the result proved in \cite{3} that states that if $ \mathbb{H} $ is a finite dimensional real Hilbert space then $ T \in \mathbb{L}(\mathbb{H}) $ is left symmetric if and only if $ T $ is the zero operator i.e., $ \mathbb{L}(\mathbb{H}) $ has no nonzero left symmetric point. \\
 
\noindent In the next theorem we prove that if $ \mathbb{X} $ is a finite dimensional strictly convex and smooth real Banach space, then a ``large" class of operators can not be left symmetric in $ \mathbb{L}(\mathbb{X}). $

\begin{theorem}
Let $ \mathbb{X} $ be a finite dimensional strictly convex and smooth real Banach space. Let $ T \in \mathbb{L}(\mathbb{X}) $ be such that there exists $ x, y \in S_{\mathbb{X}} $ satisfying $ (i)~ x \in M_T, (ii)~ y \perp_B x, (iii)~ Ty \neq 0.  $ Then $ T $ can not be left symmetric.
\end{theorem}

\noindent \textbf{Proof :} There exists a hyperplane $ H $ of codimension $ 1 $ in $ \mathbb{X} $ such that $ y \perp_B H. $ Define  a linear operator $ A \in \mathbb{L}(\mathbb{X}) $ as follows: \\
\[  Ay = Ty, A(H) = 0.  \]
Since $ y \perp_B x $ and $ \mathbb{X} $ is smooth, it follows that $ x \in H, $ i.e., $ Ax = 0. $ Since $ \mathbb{X} $ is strictly convex, as before it is easy to show that $ M_A = \{ \pm y \}. $ We observe that $ T \perp_B A, $ since $ x \in M_T ~\mbox{and}~ Tx \perp_B Ax = 0. $ However, $ M_A = \{ \pm y \}, Ay \not\perp_B Ty $ together implies that $ A \not\perp_B T. $ This completes the proof of the fact that $ T $ can not be left symmetric.

\begin{corollary}
Let $ \mathbb{X} $ be a finite dimensional strictly convex and smooth real Banach space. Let $ T \in \mathbb{L}(\mathbb{X}) $ be invertible. Then $ T $ can not be left symmetric. 
\end{corollary}
\noindent \textbf{Proof :} Since $ \mathbb{X} $ is finite dimensional, there exists $ x \in S_{\mathbb{X}} $ such that $ \| Tx \| = \| T \|. $ From Theorem $ 2.3 $ of James \cite{4}, it follows that there exists $ y (\neq 0) \in \mathbb{X} $ such that $ y \perp_{B} x. $ Using the homogeneity property of Birkhoff-James orthgonality, we may assume without loss of generality that $ \| y \| = 1. $ Since $ T $ is invertible, $ Ty \neq 0. $ Thus, all the conditions of the previous theorem are satisfied and hence $ T $ can not be left symmetric.\\
\\
In the next two theorems we establish some conditions corresponding to the right symmetry of a linear operator defined on a finite dimensional strictly convex and smooth real Banach space.

\begin{theorem}
Let $ \mathbb{X} $ be an $ n $ dimensional strictly convex and smooth real Banach space. Let $ x_0 \in S_{\mathbb{X}} $ be a left symmetric point. Let $ T \in \mathbb{L}(\mathbb{X})  $ be such that $ M_T = \{ \pm x_0 \} $ and $ x_0 $ is an eigen vector of $ T. $ Then either of the following is true: \\
\noindent (i) $ rank~ T \geq n-1. $ \\
\noindent (ii) $ T $ is not a right symmetric point in $ \mathbb{L}(\mathbb{X}). $
\end{theorem}

\noindent \textbf{Proof :} We first note that the theorem is trivially true if $ n \leq 2. $ Let $ n > 2. $ Since $ x_0 $ is an eigen vector of $ T, $ there exists a scalar $ \lambda_0 $ such that $ Tx_0 = \lambda_0 x_0. $ We also note that since $ M_T = \{ \pm x_0 \}, \lambda_0 \neq 0. $ If  $ rank~ T $ $ \geq n-1 $ then we are done. Let $ rank~ T < n-1. $ Then $ ker~ T $ is a subspace of $ \mathbb{X} $ of dimension at least $ 2. $ Let $ x_0 \perp_B H_0, $ where $ H_0 $ is a hyperplane of codimension $ 1 $ in $ \mathbb{X}. $ Since $ dim~ker~T \geq 2, $ there exists a unit vector $ u_0 \in S_{\mathbb{X}} $ such that $ u_0 \in H_0 \cap ker~T. $ Since $ x_0 $ is a left symmetric point and $ x_0 \perp_B u_0, $ we have $ u_0 \perp_B x_0. $ Since $ \mathbb{X} $ is smooth, there exists a unique hyperplane $ H_1 $ of codimension $ 1 $ in $ \mathbb{X} $ such that $ u_0 \perp_B H_1. $ Clearly, $ x_0 \in H_1. $ Let $ \{ x_0, y_i : i = 1, 2, \ldots, n-2\} $ be a  basis of $ H_1. $ Then $ \{ u_0, x_0, y_i : i = 1, 2, \ldots, n-2 \} $ is  basis of $ \mathbb{X} $ such that $ u_0 \perp_B span\{ x_0, y_i : i = 1, 2, \ldots, n-2 \}. $ Define a linear operator $ A \in \mathbb{L}(\mathbb{X}) $ as follows: \\
\[  Au_0 = u_0, Ax_0 = \frac{1}{2} x_0, Ay_i = \frac{1}{2} y_i. \]

\noindent It is routine to check that $ M_A = \{ \pm u_0 \}, $ since $ \mathbb{X} $ is strictly convex and $ u_0 \perp_B span\{ x_0, y_i : i = 1, 2, \ldots, n-2 \}. $  Since $ u_0 \perp_B x_0 $ and Birkhoff-James orthogonality is homogeneous, $ u_0 = Au_0 \perp_B \lambda_0 x_0 = Tx_0, $ which proves that $ A \perp_B T. $ However, since $ \lambda_0 \neq 0, Tx_0 = \lambda_0 x_0 \not\perp_B \frac{1}{2} x_0 = Ax_0. $ This, coupled with the fact that $ M_T = \{ \pm x_0 \}, $ implies that $ T \not\perp_B A $ and thus $ T $ is not a right symmetric point in $ \mathbb{L}(\mathbb{X}). $

\begin{theorem}
Let $ \mathbb{X} $ be an $ n $ dimensional strictly convex and smooth real Banach space. Let $ T \in \mathbb{L}(\mathbb{X}) $ be such that $ M_T = \{ \pm x_0 \} $ and $ ker~T $ contains a nonzero left symmetric point. Then either of the following is true: \\
\noindent (i) $ I \perp_B T $ and $ T \perp_B I, $ where $ I \in \mathbb{L}(\mathbb{X}) $ is the identity operator on $ \mathbb{X}. $ \\
\noindent (ii) $ T $ is not a right symmetric point in $ \mathbb{L}(\mathbb{X}). $
\end{theorem}

\noindent \textbf{Proof :} Let $ u_0 \in ker~T $ be a nonzero left symmetric point. Without loss of generality let us assume that $ \| u_0 \| = 1. $ We have, $ \| I + \lambda T \| \geq \| (I + \lambda T) u_0 \| = 1 \geq \| I \|, $ which proves that $ I \perp_B T. $ If $ T \perp_B I $ then we are done. If possible suppose that $ T \not\perp_B I. $ Since $ M_T = \{ \pm x_0 \}, $ it follows that $ Tx_0 \not\perp_B Ix_0 =x_0. $ Since $ \mathbb{X} $ is strictly convex, $ u_0 \in S_{\mathbb{X}} $ is an exposed point of the unit ball $ B_{\mathbb{X}}. $ Let $ H_0 $ be the hyperplane of codimension $ 1 $ in $ \mathbb{X} $ such that $ u_0 \perp_B H_0. $ Let $ \{u_1, u_2, \ldots, u_{n-1}\} $ be a  basis of $ H_0. $ Then $ \{u_0, u_1, \ldots, u_{n-1}\} $ is a basis of $ \mathbb{X} $ such that $ u_0 \perp_B span\{ u_1, u_2, \ldots, u_{n-1} \}. $ \\
\noindent Let $ x_0 = \alpha_0 u_0 + \alpha_1 u_1 + \ldots + \alpha_{n-1} u_{n-1}, $ for some scalars $ \alpha_0, \alpha_1, \ldots, \alpha_{n-1}.  $ Clearly, $ u_0 \perp_B{u_i} $ for each $ i = 1,2, \ldots, n-1. $ Since $ \mathbb{X} $ is smooth, Birkhoff-James orthogonality is right additive in $ \mathbb{X} $ and thus we have $ u_0 \perp_B \alpha_1 u_1 + \ldots + \alpha_{n-1} u_{n-1}. $ Since $ u_0 $ is a left symmetric point in $ \mathbb{X}, $  $ \alpha_1 u_1 + \ldots + \alpha_{n-1} u_{n-1} \perp_B u_0. $ We claim that $ \alpha_0 = 0. $ If not, then using the strict convexity of $ \mathbb{X}, $ we have, \\
\noindent $ 1 = \| x_0 \| = \| \alpha_0 u_0 + (\alpha_1 u_1 + \ldots + \alpha_{n-1} u_{n-1}) \| = \|  (\alpha_1 u_1 + \ldots + \alpha_{n-1} u_{n-1}) + \alpha_0 u_0 \| > \| (\alpha_1 u_1 + \ldots + \alpha_{n-1} u_{n-1}) \|.$ Therefore, we also have, $ \| T(\alpha_1 u_1 + \ldots + \alpha_{n-1} u_{n-1}) \| = \| T(\alpha_0 u_0 + \alpha_1 u_1 + \ldots + \alpha_{n-1} u_{n-1}) \| = \| Tx_0 \| = \| T \|,  $ a contradiction. Therefore $ \alpha_0 = 0 $ and our claim is proved. Thus we have, $ x_0 = \alpha_1 u_1 + \ldots + \alpha_{n-1} u_{n-1} $ and $ u_0 \perp_B x_0. $ \\
\noindent Since $ \mathbb{X} $ is smooth, $ H_0 $ is unique in the sense that if $ H $ is any hyperplane of codimension $ 1 $ in $ \mathbb{X} $ such that $ u_0 \perp_B H $ then $ H = H_0. $ Clearly, $ x_0 \in H_0. $ Let $ \{ u_0, x_0, y_i : i = 3, 4, \ldots , n \} $ be a basis of $ \mathbb{X} $ such that $ u_0 \perp_B span\{ x_0, y_i : i = 3, 4, \ldots , n \} $. Define a linear operator $ A \in \mathbb{L}(\mathbb{X}) $ as follows: 
\[ Au_0 = u_0, Ax_0 = \frac{1}{2} x_0, Ay_i = \frac{1}{2} y_i. \]
\noindent As before, it is easy to check that $ M_A = \{ \pm u_0 \}, $ since $ \mathbb{X} $ strictly convex and $ u_0 \perp_B span\{ x_0, y_i : i = 3, 4, \ldots , n \} $. Clearly, $ A \perp_B T, $ since $ Au_0 \perp_B Tu_0 = 0. $ We also note that since $ M_T = \{ \pm x_0 \} $ and $ Tx_0 \not\perp_B Ax_0 = \frac{1}{2} x_0, $ we must have that $ T \not\perp_B A. $ This proves that $ T $ is not a right symmetric point in $ \mathbb{L}(\mathbb{X}) $ and completes the proof of the theorem.\\

\noindent Let $ \mathbb{H} $ be a finite dimensional real Hilbert space. It was proved in \cite{3} that $ T \in \mathbb{L}(\mathbb{H}) $ is a left symmetric point in $ \mathbb{L}(\mathbb{H}) $ if and only if $ T $ is the zero operator. In the next example we show that if $ \mathbb{X} $ is a finite dimensional real Banach space, which is not a Hilbert space, then there may exist nonzero left symmetric operators in $ \mathbb{L}(\mathbb{X}). $

\begin{example}
Let $ \mathbb{X} $ be the $ 2 $ dimensional real Banach space $ l_{1}^{2}. $ Let $ T \in \mathbb{L}(\mathbb{X}) $ be defined by $ T(1, 0) = (\frac{1}{2}, \frac{1}{2}), T(0, 1) = (0, 0). $ We claim that $ T $ is left symmetric in $ \mathbb{L}(\mathbb{X}). $ Indeed, let $ A \in \mathbb{L}(\mathbb{X}) $ be such that $ T \perp_B A. $ Since $ M_T = \{ \pm (1, 0) \}, $ it follows that $ (\frac{1}{2}, \frac{1}{2}) = T(1, 0) \perp_B A(1, 0). $ Since $ (\frac{1}{2}, \frac{1}{2}) $ is a left symmetric point in $ \mathbb{X}, $ it follows that $ A(1, 0) \perp_B T(1, 0). $ We also note that $ \{ \pm (1, 0), \pm (0, 1) \} $ are the only extreme points of $ S_{\mathbb{X}}. $ Since a linear operator defined on a finite dimensional Banach space must attain norm at some extreme point of the unit sphere, either $ (1, 0) \in M_A $ or $ (0, 1) \in M_A. $ In either case, there exists a unit vector $ x $ such that $ x \in M_A $ and $ Ax \perp_B Tx. $ However, this clearly implies that $ A \perp_B T, $ completing the proof of the fact that $ T $ is a nonzero left symmetric point in $ \mathbb{L}(\mathbb{X}). $
\end{example}

\noindent We next wish to prove that in case of the strictly convex and smooth real Banach spaces $ l_{p}^2  (p \geq 2, p \neq \infty), T \in \mathbb{L}(l_{p}^2) (p \geq 2, p \neq \infty) $ is left symmetric if and only if $ T $ is the zero operator. Before proving the desired result, we first state two easy propositions. It may be noted that the proofs of both the propositions follow easily from ordinary calculus.\\
\begin{prop} Let $ \mathbb{X}  $ be the real Banach space $ l_{p}^2 (p \neq 1, \infty). $ $ x \in S_{\mathbb{X}} $ is a left symmetric point in $ \mathbb{X} $ if and only if $ x \in \pm \{ (1, 0), (0, 1), (\frac{1}{2^{1/p}}, \frac{1}{2^{1/p}}), (\frac{1}{2^{1/p}}, \frac{-1}{2^{1/p}}) \}. $\\
\\
\end{prop}
\begin{prop} Let $ \mathbb{X}  $ be the real Banach space $ l_{p}^2 (p \neq 1, \infty). $ If $ x,y \in S_{\mathbb{X}} $ are such that $ x \perp_B y $ and $ y \perp_B x $ then either of the following is true: \\ 
\noindent (i) $ x = \pm(1, 0) $ and $ y = \pm(0, 1). $ \\
\noindent (ii) $ x = \pm(0, 1) $ and $ y = \pm(1, 0). $ \\
\noindent (iii) $ x = \pm(\frac{1}{2^{1/p}}, \frac{1}{2^{1/p}}), y = \pm(\frac{1}{2^{1/p}}, \frac{-1}{2^{1/p}}). $\\
\noindent (iv) $ x = \pm(\frac{1}{2^{1/p}}, \frac{-1}{2^{1/p}}), y = \pm(\frac{1}{2^{1/p}}, \frac{1}{2^{1/p}}). $ 
\\
\end{prop}

\noindent Next we apply these two propositions and some of the results proved in this paper to prove that $ T \in \mathbb{L}(l_{p}^2) (p \geq 2, p \neq \infty) $ is left symmetric if and only if $ T $ is the zero operator.

\begin{theorem}
Let $ \mathbb{X} $ be the $ 2 $ dimensional real Banach space $ l_{p}^2 (p \geq 2, p \neq \infty) $. $ T \in \mathbb{L}(\mathbb{X}) $ is left symmetric if and only if $ T $ is the zero operator.
\end{theorem}

\noindent \textbf{Proof :} If possible suppose that $ T \in \mathbb{L}(\mathbb{X}) $ is a nonzero left symmetric point in $ \mathbb{L}(\mathbb{X}). $ Since Birkhoff-James orthogonality is homogeneous, and $ T $ is nonzero, let us assume, without loss of generality, that $ \| T \| = 1. $ Let $ T $ attains norm at $ x \in S_{\mathbb{X}}. $ From $ \mbox{Theorem} ~ 2.3 $ of James \cite{4}, it follows that there exists $ y \in S_{\mathbb{X}} $ such that $ y \perp_B x. $ Since $ \mathbb{X} $ is strictly convex and smooth, applying $ \mbox{Theorem}~ 2.5,  $ we see that $ Ty = 0. $ $ \mbox{Theorem}~ 2.4 $ ensures that $ Tx $ must be a left symmetric point in $ \mathbb{X}. $ Thus, applying $ \mbox{Proposition} ~ 2.8, $ we have that 
\[ Tx \in \pm \{ (1, 0), (0, 1), (\frac{1}{2^{1/p}}, \frac{1}{2^{1/p}}), (\frac{1}{2^{1/p}}, \frac{-1}{2^{1/p}}) \}. \]

\noindent We next claim that $ x \perp_B y. $ \\
\noindent From $ \mbox{Theorem}~ 2.3 $ of James \cite{4}, it follows that there exists a real number $ a $ such that $ ay + x \perp_B y. $ Since $ y \perp_B x $ and $ x, y \neq 0, $ $ \{ x, y \} $ is linearly independent and hence $ ay + x \neq 0. $ Let $ z = \frac{ay + x}{\| ay + x \|}. $ We note that if $ Tz = 0 $ then $ T $ is the zero operator. Let $ Tz \neq 0. $ Clearly, $ \{ y, z \} $ is a basis of $ \mathbb{X}. $ \\
\noindent Let $ \| c_1 z + c_2 y \| = 1, $ for some scalars $ c_1, c_2. $ Then we have, $ 1 = \| c_1 z + c_2 y \| \geq ~\mid c_1 \mid. $ Since $ \mathbb{X} $ is strictly convex, $ 1 > \mid c_1 \mid, $ if $ c_2 \neq 0. $ We also have, $ \| T(c_1 z + c_2 y) \| = \| c_1 Tz \| = \mid c_1 \mid \| Tz \| \leq \| Tz \| $ and $ \| T(c_1 z + c_2 y) \| = \| Tz \| $ if and only if $ c_1 = \pm 1 $ and $ c_2 = 0. $ This proves that $ M_T = \{ \pm z \}. $ However, we have already assumed that $ x \in M_T. $ Thus, we must have $ x = \pm z. $ Since $ z \perp_B y, $ our claim is proved.\\
\noindent Thus, $ x, y \in S_{\mathbb{X}} $ are such that $ x \perp_B y $ and $ y \perp_B x. $ Therefore, by applying Proposition $ 2.9, $ we see that we have the following information about $ T : $ \\ 
\noindent (i) $ T $ attains norm at $ x, $ $ x \perp_B y, y \perp_B x, Ty = 0. $ \\
\noindent (ii) $ x, y, Tx \in \pm \{ (1, 0), (0, 1), (\frac{1}{2^{1/p}}, \frac{1}{2^{1/p}}), (\frac{1}{2^{1/p}}, \frac{-1}{2^{1/p}}) \}. $ \\
This effectively ensures that in order to prove that $ T \in \mathbb{L}(\mathbb{X}) $ is left symmetric if and only if $ T $ is the zero operator, we only need to consider $ 32 $ different  operators that satisfy (i) and (ii) and show that none among them is left symmetric. \\
\noindent Let us first consider one such typical linear operator and prove that it is not left symmetric. \\ 
\noindent Let $ T \in \mathbb{L}(\mathbb{X}) $ be defined by: $ T(1, 0) = (1, 0), T(0,1) = (0, 0). $ Define $ A \in \mathbb{L}(\mathbb{X}) $ by $ A(1, 0) = (0, 1), A(0, 1) = (1, 1). $ Clearly, $ T $ attains norm only at $ \pm (1, 0). $ Since $ T(1, 0) = (1, 0) \perp_B (0, 1) = A(1, 0), $ it follows that $ T \perp_B A. $ We claim that $ A \not\perp_B T. $ \\
\noindent Now, $ \| A(\frac{1}{2^{1/p}}, \frac{1}{2^{1/p}}) \|^{p} = \frac{1}{2} + 2^{p-1} > 2 = \| A(0, 1) \|^{p}, $ since $ p \geq 2. $ This proves that $ \pm (1, 0), (0, 1) \notin M_A. $ It is also easy to observe that if $ (\alpha, \beta) \in M_A $ then either $ \alpha, \beta \geq 0 $ or $ \alpha, \beta \leq 0. $ \\
\noindent For any $ \alpha, \beta > 0 $ and for sufficiently small negative $ \lambda, $ \\
\noindent $ \| A(\alpha, \beta) + \lambda T(\alpha, \beta) \|^p = \| (\beta + \lambda \alpha, \alpha + \beta) \|^p = | \beta + \lambda \alpha |^p + | \alpha + \beta |^p < | \beta |^p + | \alpha + \beta |^p = \| A(\alpha, \beta) \|^p. $ \\
\noindent Similarly,  for any $ \alpha, \beta < 0, $ and for sufficiently small negative $ \lambda,  \| A(\alpha, \beta) + \lambda T(\alpha, \beta) \|^p < \| A(\alpha, \beta) \|^p. $ \\
\noindent This proves that for any $ w \in M_A, Tw \notin Aw^{-}. $ Applying $ \mbox{Theorem} ~ 2.2, $ it now follows that $ A \not\perp_B T. $ Thus, $ T $ is not left symmetric in $ \mathbb{L}(\mathbb{X}), $ contradicting our initial assumption. \\
\noindent Next we describe a general method to prove that none among these $ 32 $ operators are left symmetric. \\
\noindent Let $ T $ attains norm at $ x, $ $ x \perp_B y, y \perp_B x, Ty = 0 $ and $ x, y, Tx \in \pm \{ (1, 0), (0, 1), (\frac{1}{2^{1/p}}, \frac{1}{2^{1/p}}), (\frac{1}{2^{1/p}}, \frac{-1}{2^{1/p}}) \}. $ Define a linear operator $ A \in \mathbb{L}(\mathbb{X}) $ by $ Ax = y, Ay = (1, 0)~ or~ (1, 1)  $ such that \\
\noindent (i) $ A $ does not attain its norm at $ \pm x, \pm y. $ \\
\noindent (ii) $ Tw \notin Aw^{-} $ for any $ w \in M_A. $ \\
\noindent Then, as before, it is easy to see that $ T \perp_B A  $ but $ A \not\perp_B T. $ Thus, $ T $ is not left symmetric.\\
This completes the proof of the theorem.

\begin{remark}
The general method described in the previous theorem to prove that $ T \in \mathbb{L}(l_{p}^2) (p \geq 2, p \neq \infty) $ is left symmetric if and only if $ T $ is the zero operator has been verified separately in each possible case mentioned in the proof of the theorem. In this paper, only one particular case has been dealt with explicitly. The details have been omitted in other cases since the method remains the same in each case.
\end{remark}

\begin{remark}
We intuitively believe the above result to hold for any $ l_{p}^{n} (p > 1, p \neq \infty)$ space and more generally for any finite dimensional strictly convex and smooth real Banach space. In fact, motivated by the ideas used in the proof of the previous theorem, we make the following two conjectures:
\end{remark}

\begin{conj}
Let $ \mathbb{X} $ be the  finite dimensional real Banach space $ l_{p}^n (p > 1, p \neq \infty) $. $ T \in \mathbb{L}(\mathbb{X}) $ is left symmetric if and only if $ T $ is the zero operator.
\end{conj}

\begin{conj}
Let $ \mathbb{X} $ be a finite dimensional strictly convex and smooth real Banach space. $ T \in \mathbb{L}(\mathbb{X}) $ is left symmetric if and only if $ T $ is the zero operator.
\end{conj}

\noindent Of course, if Conjecture 2.14 is true then so is Conjecture 2.13. Conversely, if Conjecture 2.13 is false then so is Conjecture 2.14. Yet, we mention them separately because it seems to us that proving or disproving Conjecture 2.13 might be easier and is interesting in its own right. We would also like to add that exploring and characterizing the symmetry of Birkhoff-James orthogonality of linear operators defined on a finite dimensional strictly convex and smooth real Banach space seems to be a challenging as well as compelling problem, which is all the more relevant in light of the results obtained in the present paper.\\

\noindent \textbf{Acknowledgement.} The author would like to lovingly acknowledge the contribution of Dr. Dwijendra Nath Sain, his father and his first teacher, towards shaping his mathematical philosophy and for the constant inspiration and encouragement. He would also like to thank Dr. Kallol Paul and Ms. Puja Ghosh for their invaluable suggestions while writing this paper.
 
\bibliographystyle{amsplain}

\end{document}